\documentclass[12pt]{amsart}
\usepackage{amssymb,amsmath,amsthm,mathtools,amsfonts,times,hyperref,mathrsfs,multirow,tabularx}
\usepackage{enumerate}
\usepackage{graphicx}
\usepackage{array}
\usepackage{ytableau}
\usepackage{pstricks,pst-node,graphicx}
\usepackage{tikz,varwidth}
\usepackage{setspace}
\usepackage{float}
\usetikzlibrary{calc}
\usepackage{tikz}
\usepackage{extarrows}
\usepackage[super]{nth}

\newtheorem{theorem}{Theorem}[section]

\newtheorem{conj}[theorem]{Conjecture}

\theoremstyle{definition}

\theoremstyle{remark}
\newtheorem{remark}{Remark}
\numberwithin{equation}{section}

\allowdisplaybreaks

\newcommand\numberthis{\addtocounter{equation}{1}\tag{\theequation}}
\newcommand\restr[2]{{
  \left.\kern-\nulldelimiterspace 
  #1 
  \littletaller 
  \right|_{#2} 
  }}
\newcommand{\littletaller}{\mathchoice{\vphantom{\big|}}{}{}{}}  

\begin{document}

\title[]
{New Borwein-type conjectures}
    
\author{Alexander Berkovich}
\address{Department of Mathematics, University of Florida, Gainesville
FL 32611, USA}
\email{alexb@ufl.edu}
\author{Aritram Dhar}
\address{Department of Mathematics, University of Florida, Gainesville
FL 32611, USA}
\email{aritramdhar@ufl.edu}

\date{\today}



\begin{abstract}
Motivated by recent research of Wang and Krattenthaler, we use Maple to propose five new ``Borwein-type'' conjectures modulo $3$ and two new ``Borwein-type'' conjectures modulo $5$.
\end{abstract}
\maketitle

\section{Introduction}\label{s1}
In $1990$, Peter Borwein unearthed several mysteries involving sign patterns of coefficients of certain polynomials. He then personally communicated his observations to George Andrews who stated them as the following in his paper \cite{A95}:\\ 
\begin{conj}[First Borwein Conjecture]\label{conj11}
For the polynomials $A_n(q)$, $B_n(q)$, and $C_n(q)$ defined by\\
\begin{align*}
\prod\limits_{j = 1}^{n}(1-q^{3j-2})(1-q^{3j-1}) = A_n(q^3) - qB_n(q^3) - q^2C_n(q^3),\\    
\end{align*}
each has non-negative coefficients.\\
\end{conj}
\begin{conj}[Second Borwein Conjecture]\label{conj12}
For the polynomials $\alpha_n(q)$, $\beta_n(q)$, and $\gamma_n(q)$ defined by\\
\begin{align*}
\prod\limits_{j = 1}^{n}(1-q^{3j-2})^2(1-q^{3j-1})^2 = \alpha_n(q^3) - q\beta_n(q^3) - q^2\gamma_n(q^3),\\    
\end{align*}
each has non-negative coefficients.\\    
\end{conj}
\begin{conj}[Third Borwein Conjecture]\label{conj13}
For the polynomials $\nu_n(q)$, $\phi_n(q)$, $\chi_n(q)$, $\psi_n(q)$, and $\omega_n(q)$ defined by\\
\begin{align*}
\prod\limits_{j = 1}^{n}(1-q^{5j-4})(1-q^{5j-3})(1-q^{5j-2})(1-q^{5j-1})\\ = \nu_n(q^5) - q\phi_n(q^5) - q^2\chi_n(q^5) - q^3\psi(q^5) - q^4\omega(q^5),\\    
\end{align*}
each has non-negative coefficients.\\    
\end{conj}
Although Andrews did not prove any of the above conjectures, he derived explicit sum representations \cite[Theorem $3.1$, eq. $(3.4)$ -- $(3.6)$]{A95} of the polynomials $A_n(q)$, $B_n(q)$, and $C_n(q)$ respectively. Conjecture \ref{conj11} was recently proved analytically by Chen Wang using the saddle point method \cite{W22}. Apart from other tools, his proof relied on a theorem of Andrews \cite[Theorem $4.1$]{A95} and the following three recursive relations also given by Andrews \cite[eq. $(3.1)$ -- $(3.3)$]{A95}:\\
\begin{align*}
A_n(q) &= (1+q^{2n-1})A_{n-1}(q) + q^nB_{n-1}(q) + q^nC_{n-1}(q),\\
B_n(q) &= (1+q^{2n-1})B_{n-1}(q) + q^{n-1}A_{n-1}(q) - q^nC_{n-1}(q),\\
C_n(q) &= (1+q^{2n-1})C_{n-1}(q) + q^{n-1}A_{n-1}(q) - q^{n-1}B_{n-1}(q).\\
\end{align*}
\par Recently, Chen Wang and Christian Krattenthaler \cite{KW22} proved Conjecture \ref{conj11} and Conjecture \ref{conj12}. Unlike Wang \cite{W22} who considered Andrews' sum representations for his proof, Wang and Krattenthaler \cite{KW22} applied the saddle point method directly to the ``First Borwein polynomial'' $\prod\limits_{j = 1}^{n}(1-q^{3j-2})(1-q^{3j-1})$ and its higher powers. They also gave \cite[Conjecture $1.3$]{KW22} a cubic conjecture which was missed by both Borwein and Andrews and can be stated as follows:\\
\begin{conj}[Cubic Borwein Conjecture]\label{conj14}
For the polynomials $\kappa_n(q)$, $\delta_n(q)$, and $\theta_n(q)$ defined by\\
\begin{align*}
\prod\limits_{j = 1}^{n}(1-q^{3j-2})^3(1-q^{3j-1})^3 = \kappa_n(q^3) - q\delta_n(q^3) - q^2\theta_n(q^3),\\    
\end{align*}
each has non-negative coefficients.\\    
\end{conj}
They proved ``two thirds'' of Conjecture \ref{conj14} above by showing that all coefficients of $\kappa_n(q)$ are non-negative, ``one half'' of the coefficients of $\delta_n(q)$ are non-negative, and ``one half'' of the coefficients of $\theta_n(q)$ are non-negative. They also provided ideas \cite[Conjecture $11.1$]{KW22} along similar lines for (at least) ``three fifths'' of a proof of Conjecture \ref{conj13} when the modulus $3$ is replaced by $5$.\\\par From here on, let us consider the standard notation for $q$-Pochammer symbols,\\
\begin{align*}
(a;q)_{n} &= \prod\limits_{k = 0}^{n-1}(1-aq^k),\,\,\text{for $n\ge 1$},\\
(a;q)_{0} &= 1.\\
\end{align*}
\par For positive integers $n$ and $k$, we consider polynomials of the form\\
\begin{align*}
P_{n,k}(q) := \dfrac{(q;q)_{kn}}{(q^k;q^k)_n} = \prod\limits_{j=1}^{k-1}(q^j;q^k)_n.\\    
\end{align*}
\par For $k = 7$ above, Wang and Krattenthaler \cite[Conjecture $11.5$]{KW22} stated an interesting modulus $7$ Borwein Conjecture which is as follows:\\
\begin{conj}[Modulus $7$ Borwein Conjecture]\label{conj15}
For positive integers $n$, consider the expansion of the polynomial\\
\begin{align*}
P_{n,7}(q) := \dfrac{(q;q)_{7n}}{(q^7;q^7)_n} = \sum\limits_{m = 0}^{21n^2}d_m(n)q^m.   
\end{align*}
Then\\
\begin{align*}
d_{7m}(n)\ge 0\quad\,\text{and}\quad\,d_{7m+1}(n),d_{7m+3}(n),d_{7m+4}(n),d_{7m+6}(n)\le 0,\quad\,\text{for all $m$ and $n$},\\    
\end{align*}
while\\
$$d_{7m+5}(n) \begin{cases}
\ge 0,\quad\text{for}\,\,\,m\le3\alpha(n)n^2,\\
\le 0,\quad\text{for}\,\,\,m > 3\alpha(n)n^2,\\
\end{cases}$$\\
where $\alpha(n)$ seems to stabilise around $0.302$.\\
\end{conj}
\par Motivated by Conjecture \ref{conj15}, in Section \ref{s2}, we state new ``Borwein-type'' sign-pattern conjectures for the polynomial $P_{n,k}^{\,i}(q)$ where $(k,i)\in\{(3,4),(3,5),(3,6),(3,7),(3,8),(5,2),(5,3)\}$. We conclude with some possible future directions in Section \ref{s3}.\\

\section{New Conjectures}\label{s2}
In this section, we state the new ``Borwein-type'' conjectures.\\
\begin{conj}\label{conj21}
Let $n$ be a positive integer and $i\in\{4,5,6,7,8\}$. Consider the expansion of the polynomial $P_{n,3}^{\,i}(q)$ defined by\\
\begin{align*}
P_{n,3}^{\,i}(q) := \dfrac{(q;q)_{3n}^i}{(q^3;q^3)_n^i} = \sum\limits_{m = 0}^{3in^2}c_m^{(i)}(n)q^m.    
\end{align*}
Then\\
\begin{align*}
c_{3m}^{(i)}(n)\ge 0,\quad\,\text{for all $m$ and $n$},    
\end{align*}
while\\
$$c_{3m+2}^{(i)}(n) \begin{cases}
\ge 0,\quad\text{for}\,\,\,m\le\alpha^{(i)}(n)n^2,\\
\le 0,\quad\text{for}\,\,\,m > \alpha^{(i)}(n)n^2,\\
\end{cases}$$\\
where $\alpha^{(i)}(n)$ tends to a limiting value as $n$ tends to $\infty$ according to the following table\\
\begin{center}
\begin{tabularx}{1.0\textwidth} { 
| >{\centering\arraybackslash}X 
| >{\centering\arraybackslash}X | }
\hline
value of $i$ & $\lim\limits_{n\rightarrow\infty}\alpha^{(i)}(n)$\\
\hline
$4$ & $0.750\cdots$ \\
\hline
$5$ & $1.270\cdots$ \\
\hline
$6$ & $1.778\cdots$ \\
\hline
$7$ & $2.282\cdots$ \\
\hline
$8$ & $2.784\cdots$ \\
\hline
\end{tabularx}\\\quad\\\quad\\
\end{center}
\end{conj}
We verified Conjecture \ref{conj21}\,\, using $q$-series package in Maple. For instance, we got\\
\begin{align*}
    \alpha^{(4)}(65)\cdot{65^2} = 3170.
\end{align*}\\
\vspace{-8mm}
\begin{remark}\label{rmk1}
Since the polynomial $P_{n,3}^{\,i}(q)$ is palindromic for any positive integer $i$, we can also predict the signs of the coefficients $c_{3m+1}^{(i)}(n)$ in Conjecture \ref{conj21}. This is true because of the following duality result.\\
\end{remark}
Consider the polynomial $P_{n,3}^{\,i}(q)$ as defined in Conjecture \ref{conj21}. If we write\\
\begin{align*}
  P_{n,3}^{\,i}(q) = X_{n,i}(q^3) + qY_{n,i}(q^3) + q^2Z_{n,i}(q^3),\numberthis\label{eq21}
\end{align*}\\
then it is easy to show that the following relations hold\\
\begin{align*}
  X_{n,i}(q) = q^{in^2}X_{n,i}(1/q)\numberthis\label{eq22}
\end{align*}\\
and\\
\begin{align*}
  Z_{n,i}(q) = q^{in^2-1}Y_{n,i}(1/q).\numberthis\label{eq23}  
\end{align*}\\
\begin{conj}\label{conj22}
Let $n$ be a positive integer and $i\in\{2,3\}$. Consider the expansion of the polynomial $P_{n,5}^{\,i}(q)$ defined by\\
\begin{align*}
P_{n,5}^{\,i}(q) := \dfrac{(q;q)_{5n}^i}{(q^5;q^5)_n^i} = \sum\limits_{m = 0}^{10in^2}d_m^{(i)}(n)q^m.    
\end{align*}
Then\\
\begin{align*}
d_{5m}^{(i)}(n)\ge 0,\quad\,\text{for all $m$ and $n$},    
\end{align*}
while\\
$$d_{5m+3}^{(i)}(n) \begin{cases}
\ge 0,\quad\text{for}\,\,\,m\le\Tilde{\alpha}^{(i)}(n)n^2,\\
\le 0,\quad\text{for}\,\,\,m > \Tilde{\alpha}^{(i)}(n)n^2,\\
\end{cases}$$\\
and\\
$$d_{5m+4}^{(i)}(n) \begin{cases}
\ge 0,\quad\text{for}\,\,\,m\le\Tilde{\beta}^{(i)}(n)n^2,\\
\le 0,\quad\text{for}\,\,\,m > \Tilde{\beta}^{(i)}(n)n^2,\\
\end{cases}$$\\
where $\Tilde{\alpha}^{(i)}(n)$ and $\Tilde{\beta}^{(i)}(n)$ tend to limiting values as $n$ tends to $\infty$ according to the following table\\
\begin{center}
\begin{tabularx}{1.0\textwidth} { 
| >{\centering\arraybackslash}X 
| >{\centering\arraybackslash}X
| >{\centering\arraybackslash}X | }
\hline
value of $i$ & $\lim\limits_{n\rightarrow\infty}\Tilde{\alpha}^{(i)}(n)$ & $\lim\limits_{n\rightarrow\infty}\Tilde{\beta}^{(i)}(n)$\\
\hline
$2$ & $1.133\cdots$ & $1.011\cdots$ \\
\hline
$3$ & $2.132\cdots$ & $2.001\cdots$ \\
\hline
\end{tabularx}\\\quad\\\quad\\
\end{center}
\end{conj}
We verified Conjecture \ref{conj22}\,\, using $q$-series package in Maple. For instance, we got\\
\begin{align*}
    \Tilde{\beta}^{(2)}(65)\cdot{65^2} = 4274.
\end{align*}\\
\vspace{-10mm}
\begin{remark}\label{rmk2}
Since the polynomial $P_{n,5}^{\,i}(q)$ is palindromic for any positive integer $i$, we can also predict the signs of the coefficients $d_{5m+1}^{(i)}(n)$ and $d_{5m+2}^{(i)}(n)$ in Conjecture \ref{conj22}. This is true because of the following duality result.\\
\end{remark}
Consider the polynomial $P_{n,5}^{\,i}(q)$ as defined in Conjecture \ref{conj22}. If we write\\
\begin{align*}
    P_{n,5}^{\,i}(q) = A_{n,i}(q^5) + qB_{n,i}(q^5) + q^2C_{n,i}(q^5) + q^3D_{n,i}(q^5) + q^4E_{n,i}(q^5),\numberthis\label{eq24}  
\end{align*}\\
then it is easy to show that the following relations hold\\
\begin{align*}
    A_{n,i}(q) = q^{2in^2}A_{n,i}(1/q),\numberthis\label{eq25}
\end{align*}\\
\begin{align*}
    D_{n,i}(q) = q^{2in^2-1}C_{n,i}(1/q),\numberthis\label{eq26}
\end{align*}\\
and\\
\begin{align*}
    E_{n,i}(q) = q^{2in^2-1}B_{n,i}(1/q).\numberthis\label{eq27}
\end{align*}\\
\begin{remark}\label{rmk3}
It must be noted that Wang and Krattenthaler \cite[Remark $11.6. (2)$]{KW22} mention that for $\delta\in\{2,3\}$, the coefficient of $q^{5m}$ in $(q;q)_{5n}^\delta/(q^5;q^5)_n^\delta$ is non-negative for all $m$, however, they do not state the limiting values of $\Tilde{\alpha}^{(\delta)}(n)$ and $\Tilde{\beta}^{(\delta)}(n)$.\\
\end{remark}

\section{Concluding Remarks}\label{s3}
\begin{enumerate}
\item It would be very interesting to investigate whether the asymptotic approach of Wang and Krattenthaler in \cite{KW22} applies to Conjecture \ref{conj21} and Conjecture \ref{conj22}.\\
\item The sign-pattern of coefficients of $P_{n,3}^{\,i}(q)$ is much more complicated when $i\ge 9$. For instance, we have\\
\begin{align*}
    X_{11,9}(q) &= 1 - 3q + 8q^3 - 9q^4 + 17q^6 - 27q^7 + 46q^9 - 57q^{10} + 260q^{12}\\ &+ 1899q^{13} + \cdots + 1899q^{1076} + 260q^{1077} - 57q^{1079} + 46q^{1080} - 27q^{1082}\\ &+ 17q^{1083} - 9q^{1085} + 8q^{1086} - 3q^{1088} + q^{1089},
\end{align*}\\
where $X_{11,9}(q)$ is as defined in \eqref{eq21} with $n = 11$ and $i = 9$. Definitely the sign-pattern above is something worth exploring further.\\
\item Similarly, the sign-pattern of coefficients of $P_{n,5}^{\,i}(q)$ is complicated when $i\ge 4$. For instance, we have\\
\begin{align*}
    \qquad\, B_{2,4}(q) &= -4 - 26q - 108q^2 - 326q^3 - 748q^4 - 1536q^5 - 2644q^6 - 4352q^7 - 6200q^8\\ &- 8418q^9 - 10544q^{10} - 11369q^{11} - 12772q^{12} - 11134q^{13} - 10072q^{14}\\ &- 7762q^{15} - 3696q^{16} - 2348q^{17} + 1420q^{18} + 2281q^{19} + 2924q^{20} + 3422q^{21}\\ &+ 2300q^{22} + 1983q^{23} + 1044q^{24} + 522q^{25} + 168q^{26} - 15q^{27} - 48q^{28} - 42q^{29}\\ &- 20q^{30} - 5q^{31},
\end{align*}\\
where $B_{2,4}(q)$ is as defined in \eqref{eq24} with $n = 2$ and $i = 4$. Definitely the sign pattern above is again worth exploring more.\\
\item It is a challenge to find partition theoretic explanation of why $X_{n,i}(q)$ in \eqref{eq21} with $1\le i\le 8$ has non-negative coefficients for all $n$. Such interpretation, if any, would lead to a great insight into positivity problems in the theory of partitions and $q$-series.\\ 
\end{enumerate}

\section{Acknowledgments}\label{s4}
The authors would like to thank Christian Krattenthaler for his wonderful talk titled ``Proofs of Borwein Conjectures'' on June $7$, $2024$ at the Legacy of Ramanujan $2024$ conference to commemorate the \nth{85} birthdays of George E. Andrews and Bruce C. Berndt held at Penn State University from June $6$-$9$, $2024$ which inspired them to pursue this topic. The authors would also like to thank the organizers of this conference: Amita Malik, James Sellers, Andrew Sills, and Ae Ja Yee for their kind invitation.\\

\bibliographystyle{amsplain}


\end{document}